\documentclass{article}
\usepackage[latin1]{inputenc}
\usepackage[T1]{fontenc}
\usepackage{amstext}
\usepackage{amsmath, amssymb}
\usepackage[english]{babel}
\usepackage{amsfonts}

\newcommand{\R}{\mathbb{R}}
\newcommand{\C}{\mathbb{C}}
\newcommand{\D}{\mathbb{D}}
\newcommand{\N}{\mathbb{N}}
\newcommand{\E}{\mathbb{E}}
\newcommand{\Z}{\mathbb{Z}}
\newcommand{\Q}{\mathbb{Q}}
\newcommand{\pp}{\mathbb{P}}
\newcommand\kD{\mathcal{D}}
\newcommand\kC{\mathcal{C}}
\newcommand\kR{\mathcal{R}}

\newcommand\kF{\mathcal{F}}
\newcommand\kE{\mathcal{E}}

\newcommand\kS{\mathcal{S}}

\newcommand\kL{\mathcal{L}}
\newcommand\ka{\mathfrak{a}}
\newcommand\kh{\mathfrak{h}}

\newtheorem {lem} {Lemma} [section]
\newtheorem {prop} {Proposition} [section]

\newtheorem {theo} {Theorem} [section]
\newtheorem {cor} {Corollary} [section]
\newtheorem {rem} {Remark} [section]

\newcommand\la{\lambda}

\title{The Heckman-Opdam Markov processes}
\author{Bruno Schapira\footnote{partially supported by the European Commission (IHP Network HARP $2002-2006$).}}

\begin{document}
\maketitle

\begin{center}
\it{Université d'Orléans,\\ Fédération Denis Poisson, Laboratoire MAPMO \\
B.P. 6759, 45067 Orléans cedex 2, France .}
\end{center}
\begin{center}
\it{Université Pierre et Marie Curie,\\
Laboratoire de Probabilités et Modèles Aléatoires, \\
4 place Jussieu, \\
F-75252 Paris cedex 05, France.}
\end{center}

\vspace*{0.6cm}

\begin{abstract}
We introduce and study the natural counterpart of the Dunkl Markov
processes in a negatively curved setting. We give a semimartingale
decomposition of the radial part, and some properties of the
jumps. We prove also a law of large numbers, a central limit
theorem, and the convergence of the normalized process to the
Dunkl process. Eventually we describe the asymptotic behavior of
the infinite loop as it was done by Anker, Bougerol and Jeulin in
the symmetric spaces setting in \cite{ABJ}.
\end{abstract}


\bigskip
\noindent{\bf A.M.S. Classification:} 58J65, 60B15, 60F05, 60F17,
60J35, 60J60, 60J65, 60J75.

\bigskip
\noindent{\bf e-mail:} bruno.schapira@univ-orleans.fr

\vspace*{0.8cm}

\section{Introduction}

In the last few years, some processes living in cones have played
an important role in probability. The cones that we consider here
are associated to a root system. Roughly speaking a root system is
a set of vectors, satisfying a few conditions, in a Euclidean
space. The set of hyperplanes orthogonal to the vectors of the
root system delimit cones, which are called the Weyl chambers. We
usually choose arbitrarily one of them which we call the positive
Weyl chamber. One of the first example of process with value in a
Weyl chamber is the intrinsic Brownian motion introduced by Biane
in \cite{Bia}. It may be defined as the radial part (in the Lie
algebras terminology) of the Brownian motion on a complex
Riemannian flat symmetric space. In the particular case where the
root system is of type $A_n$, it is the process of eigenvalues of
the Brownian motion on Hermitian matrices with trace null. It was
also proved recently by Biane, Bougerol and O'Connell \cite{BBO},
that it is a natural generalization of the Bessel-$3$ process in
dimension $n$, in the sense that it may be obtained by a transform
of the Brownian motion in $\R^n$, which coincides in dimension $1$
with the Pitman transform $2S-B$. Another example of processes
associated to a root system is the radial part of the Brownian
motion on a Riemannian symmetric space of noncompact type
\cite{ABJ}, \cite{Bab}. This is the analogue of the intrinsic
Brownian motion in a negatively curved setting. In \cite{ABJ}
Anker, Bougerol and Jeulin study in fact other processes naturally
attached to a symmetric space, and they show some surprising link
between them and the intrinsic Brownian motion. More recently,
Rösler \cite{Ros} and Rösler, Voit \cite{RV} have introduced a new
type of processes related to Weyl chambers of root systems, the
Dunkl processes. These processes are Markov processes as well as
martingales, but they no longer have continuous paths. They may
jump from a chamber to another. Nevertheless the projection on the
positive Weyl chamber of these processes, which is called the
radial part, has continuous paths. Moreover for a particular
choice of the parameter, this radial part is in fact the intrinsic
Brownian motion. These processes were studied recently more deeply
by Gallardo and Yor \cite{GY} \cite{GY2} \cite{GY3}, and by
Chybiryakov \cite{Chy} \cite{Chy2}, who have obtained many
interesting properties, such as the time inversion property, a
Wiener chaos decomposition, or a skew product decomposition.
\newline In this paper, we introduce and study the natural
counterpart of the Dunkl processes in the negatively curved
setting, which we call the Heckman-Opdam processes. These
processes are also discontinuous, and have a continuous radial
part, which coincides with the radial part (in the Lie groups
terminology) of the Brownian motion on some symmetric spaces for
particular choices of the parameter. We will show that many known
results in probability theory (see \cite{ABJ} \cite{Bab}
\cite{BJ}) in the symmetric spaces setup, can be generalized to
these new processes. We also prove that the Dunkl processes are
limits of normalized Heckman-Opdam processes.

\section{Preliminaries}

Let $\ka$ be a Euclidean vector space of dimension $n$, equipped
with an inner product $(\cdot ,\cdot )$. Let $\kh=\ka \otimes_\R
\C$ be the complexification of $\ka$. For $\alpha \in \ka$ let
$\alpha^\vee=\frac{2}{|\alpha|^2}\alpha$, and let
$$r_\alpha(x)=x-(\alpha^\vee,x)\alpha,$$
be the corresponding orthogonal reflection. Let $\kR\subset \ka$
be an integral (or crystallographic) root system, which by
definition (cf \cite{Bou}) satisfies the following hypothesis
\begin{enumerate}
\item $\kR$ is finite, does not contain $0$ and
generates $\ka$.
\item $\forall \alpha \in \kR$, $r_\alpha (\kR)=\kR$.
\item $\forall \alpha \in \kR$, $\alpha^\vee(\kR)\subset \Z$.
\end{enumerate}
We choose a set of positive roots $\kR^+$ (it can be taken as the
subset of roots $\alpha\in\kR$ such that $(\alpha,u)>0$, for some
arbitrarily chosen vector $u\in \ka$ satisfying $(\alpha,u)\neq
0$, for all $\alpha \in \kR$). We denote by $W$ the Weyl group
associated to $\kR$, i.e. the group generated by the $r_\alpha$'s,
with $\alpha \in \kR$. If $C$ is a subset of $\ka$, we call
\emph{symmetric} of $C$ any image of $C$ under the action of $W$.
Let $k\ :\ \kR \rightarrow \R^+$ be a multiplicity function, i.e.
a $W$-invariant function on $\kR$. We will assume in this paper
that $k(\alpha)$ (also denoted by $k_\alpha$ in the sequel) is
strictly positive for all $\alpha\in \kR^+$.
\newline
Let
$$\ka_+ = \{x\in \ka  \mid \forall \alpha \in \kR^+,\ (\alpha,x)>0\},$$
be the positive Weyl chamber. We denote by $\overline{\ka_+}$ its
closure, and by $\partial \ka_+$ its boundary. Let also
$\ka_{\text{reg}}$ be the subset of regular elements in $\ka$,
i.e. those elements which belong to no hyperplane $\{x\in \ka\mid
(\alpha,x)=0\}$.
\newline
For $\xi \in \ka$, let $T_\xi$ be the Dunkl-Cherednik operator. It
is defined, for $f\in C^1(\ka)$ and $x\in \ka_{\text{reg}}$, by
$$T_\xi f(x)=\partial_\xi f(x)-(\rho,\xi)f(x) + \sum_{\alpha \in \kR^+}k_\alpha
\frac{(\alpha,\xi)}{1-e^{-(\alpha,x)}}\{f(x)-f(r_\alpha x) \},$$
where
$$\rho=\frac{1}{2}\sum_{\alpha \in \kR^+}k_\alpha \alpha.$$
The Dunkl-Cherednik operators form a commutative family of
differential-difference operators (see \cite{C}). The Laplacian
$\kL$ is defined by
$$\kL=\sum_{i=1}^{n} T_{\xi_i}^2,$$
where $\{\xi_1,\dots,\xi_n\}$ is any orthonormal basis of $\ka$
($\kL$ is independent of the chosen basis). Here is an explicit
expression of $\kL$ (see \cite{Schap}), which holds for $f\in
C^2(\ka)$ and $x\in \ka_{\text{reg}}$:
\begin{eqnarray}
\label{explicit} \kL f(x) &=& \Delta f(x)+ \sum_{\alpha \in
\kR^+}k_\alpha \coth \frac{(\alpha,x)}{2}\partial_\alpha f(x)
+|\rho|^2f(x) \\
 &+& \nonumber \sum_{\alpha \in \kR^+}k_\alpha \frac{|\alpha|^2}{4\sinh^2 \frac{(\alpha,x)}{2}}
\{f(r_\alpha x)-f(x)\},
\end{eqnarray}
where $\Delta$ denotes the Euclidean Laplacian. Let $\mu$ be the
measure on $\ka$ given by
$$d\mu(x)=\delta(x)dx,$$
where $$\delta(x)=\prod_{\alpha \in \kR^+}
|\sinh\frac{(\alpha,x)}{2}|^{2k_\alpha}.$$ Let $\la \in \kh$. We
denote by $G_\la$ the unique analytic function on $\ka$, which
satisfies the differential and difference equations
$$T_\xi G_\la=(\la,\xi)G_\la,\quad \forall \xi \in \ka$$
and which is normalized by $G_\la(0)=1$ (see \cite{O}). Let
$F_\la$ be the function defined for $x\in \ka$ by
$$F_\la(x)=\frac{1}{|W|}\sum_{w\in W} G_\la(wx).$$
These functions were introduced in \cite{HO}. Let $C_0(\ka)$ be
the space of continuous functions on $\ka$ which vanish at
infinity, and let $C_0^2(\ka)$ be its subset of twice
differentiable functions (with analogue definitions for
$\overline{\ka_+}$ in place of $\ka$). We denote by $\kC(\ka)$ the
Schwartz space on $\ka$ associated to the measure $\mu$, i.e. the
space of infinitely differentiable functions $f$ on $\ka$ such
that for any polynomial $p$, and any $N\in \N$,
$$\sup_{x\in \ka} (1+|x|)^Ne^{(\rho,x^+)}|p(\frac{\partial}{\partial x}) f(x)| <
+\infty,$$ where $x^+$ is the unique symmetric of $x$ which lies
in $\overline{\ka_+}$. We denote by $\kC(\ka)^W$ the subspace of
$W$-invariant functions, which we identify with their restriction
to $\overline{\ka_+}$. We have seen in \cite{Schap} that
$$\kD:=\frac{1}{2}(\kL-|\rho|^2),$$ densely defined on $\kC(\ka)$, has a closure on
$C_0(\ka)$, which generates a Feller semigroup $(P_t,t\ge 0)$. We
have also obtained the formula for $f\in C_0(\ka)$ and $x\in \ka$:
$$P_tf(x)= \int_{\ka}p_t(x,y)f(y)d\mu(y),$$
where $p_t(\cdot,\cdot)$ is the heat kernel. It is defined for
$x,y \in \ka$ and $t>0$, by
$$p_t(x,y)=\int_{i\ka}e^{-\frac{t}{2}(|\la|^2+|\rho|^2)}G_\la(x)G_\la(-y)d\nu(\la),$$
where $\nu$ is the asymmetric Plancherel measure (see
\cite{Schap}). We denote by $D$ the differential part of $\kD$,
which is equal for $f\in \kC(\ka)^W$ and $x\in \ka_{\text{reg}}$
to
$$Df(x)=\frac{1}{2}\Delta f(x) + (\nabla \log \delta^{\frac{1}{2}},\nabla f)(x).$$ It
has also a closure on $C_0(\overline{\ka_+})$, which generates a
Feller semigroup $(P_t^W,t\ge 0)$. It is associated to a kernel
$p_t^W$, which is defined for $x,y\in \overline{\ka_+}$ and $t>0$,
by
$$p_t^W(x,y)=\int_{i\ka}e^{-\frac{t}{2}(|\la|^2+|\rho|^2)}F_\la(x)F_\la(-y)d\nu'(\la),$$
where $\nu'$ is the symmetric Plancherel measure (see
\cite{Schap}).

\section{Definition and first properties}
\subsection{The Heckman-Opdam processes}
The Heckman-Opdam process (also denoted by HO-process) is defined
as the càdlàg Feller process $(X_t,t\ge 0)$ on $\ka$ with
semigroup $(P_t,t\ge 0)$. Remember that it is characterized as the
unique (in law) solution of the martingale problem associated to
$(\kD,\kC(\ka))$ on $C_0(\ka)$, see e.g. Theorem $4.1$ and
Corollary $4.3$ in \cite{EK}. Observe that, by elementary
calculation, the generator of the HO-process is also the closure
of $\kD$ on $C_c^\infty(\ka)$, the space of infinitely
differentiable functions with compact support on $\ka$. The
multiplicity $k$ is called the parameter of the HO-process.
Moreover a.s., for any $t\ge 0$, $X_t \in \ka$, i.e. the exploding
time of $(X_t,t\ge 0)$ is almost surely infinite. This results for
example from Proposition $2.4$ in \cite{EK}. Similarly, we define
the radial process (or radial part of the HO-process), as the
Feller process on $\overline{\ka_+}$ with semigroup $(P_t^W,t\ge
0)$. It is also characterized as the unique solution of the
martingale problem associated to $(D,\kC(\ka)^W)$ on
$C_0(\overline{\ka_+})$. Consider now the process $(X^W_t,t\ge 0)$
on $\overline{\ka_+}$, defined as the projection on
$\overline{\ka_+}$ under the Weyl group $W$ (for any $t$, $X^W_t$
is the unique symmetric of $X_t$ which lies in
$\overline{\ka_+}$).
\begin{prop}
\label{premiereprop} The process $(X^W_t,t\ge 0)$ is the radial
process.
\end{prop}
\textbf{Proof of the proposition:} Remember that if $f\in
\kC(\ka)^W$, then $\kD f=Df \in \kC(\ka)^W$. Thus for $f\in
\kC(\ka)^W$, and $t\ge 0$,
$$f(X_t^W)-f(X_0^W)-\int_0^t Df(X_s^W)ds=f(X_t)-f(X_0)-\int_0^t\kD f(X_s)ds.$$
Therefore $(f(X_t^W)-f(X_0^W)-\int_0^t Df(X_s^W)ds,t\ge 0)$ is a
local martingale. And we conclude by the uniqueness of the
martingale problem associated to $D$. \hfill $\square$
\newline
\newline
Let us note eventually that when $2k$ equals the multiplicity
associated with a Riemannian symmetric space of noncompact type
$G/K$, then the radial HO-process coincide with the radial part of
the Brownian motion on this symmetric space (see \cite{HO}, or
\cite{Schap}). For instance if $G=SL_n$, then $k=\frac{1}{2}$ in
the real case, $k=1$ in the complex case, and $k=2$ in the
quaternionic case. So most of our results about radial processes
generalize known results of the probabilistic theory on symmetric
spaces.

\subsection{The Dunkl processes}
\label{definitiondunkl} We recall now the definition of the Dunkl
process and of its radial part. Let $\kR'$ be a reduced root
system (i.e. such that $\forall \alpha \in \kR'$, $2\alpha \notin
\kR'$), but non necessarily integral (i.e. we do not assume
condition $3$ in the definition), and let $k'$ be a multiplicity
function on $\kR'$. The Dunkl Laplacian $\kL'$ is defined for
$f\in C^2(\ka)$, and $x\in \ka_{\text{reg}}$ by
\begin{eqnarray*}
\kL' f(x)&=& \frac{1}{2}\Delta f(x)+ \sum_{\alpha \in
\kR'^+}k'_\alpha
\frac{1}{(\alpha,x)}\partial_\alpha f(x) \\
 &+& \sum_{\alpha \in \kR'^+}k'_\alpha \frac{1}{(\alpha,x)^2}\{f(r_\alpha
 x)-f(x)\}.
\end{eqnarray*}
It was proved by Rösler in \cite{Ros} that $\kL'$ defined on
$\kS(\ka)$, the classical Schwartz space on $\ka$,  is a closable
operator, which generates a Feller semigroup on $C_0(\ka)$.
Naturally the Dunkl process is the Feller process defined by this
semigroup. Now our proof (in \cite{Schap}) that the operator $D$
defined on $\kC(\ka)^W$ is closable, also holds in the Dunkl
setting. Thus we may define the radial Dunkl process as the Feller
process on $\overline{\ka_+}$ with generator the closure of
$(L',\kS(\ka)^W)$, where $L'$ is the differential part of $\kL'$,
and $\kS(\ka)^W$ is the subspace of $\kS(\ka)$ of $W$-invariant
function (identified with their restriction to
$\overline{\ka_+}$). In this way we get a new characterization of
the radial Dunkl process as solution of a martingale problem.
Naturally Proposition \ref{premiereprop} holds as well in the
Dunkl setting, thus our definition of the radial Dunkl process
agrees with the usual one. Eventually the intrinsic Brownian
motion is by definition the radial Dunkl process of parameter
$k'=1$.

\section{The radial HO-process as a Dirichlet process}
The goal of this section is to obtain an explicit semimartingale
decomposition of the radial HO-process. In the case of root
systems of type $A$, it was obtained by Cépa and Lépingle (see
\cite{Cep} and \cite{CepLep2} Theorem $2.2$). We present here
another approach, which is based on the theory of Dirichlet
processes. Our reference for this theory will be \cite{FukOT}.
\newline
We consider $(D,\kC(\ka)^W)$ as a symmetric operator on
$L^2(\overline{\ka_+},\mu)$ (simply denoted by $L^2$ in the
sequel). We have seen in \cite{Schap} that this operator is
closable. We denote by $(D,\kD_2(D))$ its closure. Its associated
semigroup is just the extension of $(P_t^W,t\ge 0)$ on $L^2$. It
is defined for $f,g\in L^2$ and $t\ge 0$ by
$$P_t^Wf(x)= \int_{\overline{\ka_+}}p_t^W(x,y)f(y)d\mu(y).$$
We denote by $\kE$ the associated Dirichlet form, and by $\kF$ its
domain ($\kD_2(D)\subset \kF$). It is determined for $f,g \in
\kD_2(D)$ by
$$\kE(f,g):=-\int_{\overline{\ka_+}}f(x)Dg(x)d\mu(x)=-\int_{\overline{\ka_+}}Df(x)g(x)d\mu(x).$$
The fact that $\kE$ is a regular Dirichlet form with special
standard core the algebra $\kC(\ka)^W$, results from the density
of this space in $C_0(\overline{\ka_+})$ (see Lemma $5.1$ in
\cite{Schap}) and Theorem $3.1.2$ in \cite{FukOT}. We have seen in
\cite{Schap} that when $f\in L^1(\ka,\mu)$, then $Gf:x \mapsto
\int_0^\infty P_t^Wf(x)dt$ is a.e. finite. In the terminology of
\cite{FukOT}, this means that the Dirichlet form $\kE$ (or the
semigroup $P_t^W$) is transient. This implies in particular that
the process tends to infinity when $t\to \infty$. In the sequel we
will prove a law of large numbers which makes this fact precise.
It implies also that we may consider $\kF$, equipped with its
inner product $\kE$, as a Hilbert space (see \cite{FukOT} chapter
$2$). For $i=1,\dots,n$, we denote by $\varphi_i :x\mapsto x_i$
the coordinate functions on $\overline{\ka_+}$. For $A>0$, let
$\varphi_i^A \in C^\infty(\overline{\ka_+})$ be a function which
coincides with $\varphi_i$ on $\{|x|\le A\}$, and which is null on
$\{|x|\ge A+1\}$.
\begin{lem}
For all  $A >0$, $\varphi_i^A \in \kF$, and for all $v\in
\kC(\ka)^W$,
$$\kE(\varphi_i^A,v)=-\int_{\overline{\ka_+}}D\varphi_i^A vd\mu.$$
\end{lem}
\textbf{Proof of the lemma:} Let $(u_n)_n\in \kC(\ka)^W$ which
converges uniformly to $\varphi_i^A$ as in Lemma $5.1$ in
\cite{Schap}. We will assume that $|u_n-\varphi_i^A|_\infty \le
\frac{1}{n}$ for all $n$. We have to prove that it is an
$\kE$-Cauchy sequence. Let $n<m$ be two integers. We have
\begin{eqnarray*} \kE(u_n-u_m,u_n-u_m) &=& \int_{d(x,\partial
\ka_+)\le \frac{1}{m}}D(u_n-u_m)(u_n-u_m)d\mu \\
&+& \int_{\frac{1}{m}\le d(x,\partial \ka_+)\le \frac{1}{n}}
D(u_n-u_m)(u_n-u_m)d\mu \\ &+& \int_{d(x,\partial \ka_+)\ge
\frac{1}{n}}D(u_n-u_m)(u_n-u_m)d\mu.
\end{eqnarray*}
By Lemma $5.1$ in \cite{Schap}, the integrand in the first
integral is bounded, up to a constant, by $\frac{m}{n}$. But
$\mu(\{d(x,\partial \ka_+)\le \frac{1}{m}, |x|\le A+1\})$ is
bounded, up to a constant, by $\frac{1}{m}$. Thus the first
integral tends to $0$ when $n\to \infty$. The same argument
applies for the second integral. The third integral is naturally
bounded, up to a constant, by $\frac{1}{n}$. Now always by Lemma
$5.1$ in \cite{Schap}, the sequence $(D u_n)_n$ is dominated, up
to a constant, by $x\mapsto \frac{1}{d(x,\partial
\ka_+)}+\sum_\alpha \coth \frac{\alpha}{2}$ on $|x|\le A+1$, which
is $\mu$-integrable since $k$ is strictly positive. Thus the last
statement of the lemma is a consequence of the dominated
convergence theorem.  \hfill $\square$
\newline
\newline
The lemma implies (in the terminology of \cite{FukOT}), that the
functions $\varphi_i$ are in $\kF_{b,\text{loc}}$. Thanks to
Theorem $5.5.1$ p.$228$ of \cite{FukOT}, there exist martingale
additive functionals locally of finite energy $M^i$, and additive
functionals locally of zero energy $N^i$ such that, for every $i$,
\begin{eqnarray} \label{edsradiali} \varphi_i(X^W)=M^i+N^i.
\end{eqnarray}
For $A>0$, we denote by $\nu_i^A$ the measure defined by
$d\nu_i^A=-D\varphi_i^A d\mu$. Observe that $1_{(|x|\le
A)}d\nu_i^A(x)=\sum_\alpha k_\alpha
\coth(\frac{\alpha}{2},x)\varphi_i(\alpha) 1_{(|x|\le A)}d\mu(x)$.
We denote by $\nu_i$ the Radon measure defined by
$d\nu_i(x)=\sum_\alpha k_\alpha
\coth(\frac{\alpha}{2},x)\varphi_i(\alpha)d\mu(x)$. Thanks to
Theorem $5.5.4$ p.$229$ of \cite{FukOT}, and the preceding lemma,
we see that $N^i$ is an additive functional of bounded variation,
and that it is the unique continuous additive functional
associated to the measure $\nu^i$ ($\nu^i$ is called the Revuz
measure of $N^i$). Moreover from Theorem $5.1.3$ $(iii)$, we get
for all $i$ and all $t\ge 0$,
$$N^i_t=\sum_{\alpha \in \kR^+} k_\alpha \phi_i(\alpha)\int_0^t\coth\frac{(\alpha,X_s^W)}{2}ds.$$ In the same way, it is
immediate from Theorem $5.5.2$ p.$229$, and the identity $3.2.14$
p.$110$, that the Revuz measure of $<M^i>$ is $\mu$ for each $i$.
Therefore $M:=\sum_i M^i e_i$ is necessarily a Brownian motion on
$\ka$ (we have denoted by $e_i$ the $i^{\text{th}}$ vector of the
canonical basis). Now, since $\coth$ is positive on $(0,+\infty)$,
and $X^W$ does not explode, we see with (\ref{edsradiali}) that
necessarily, for all $t>0$, $\sum_{i=1}^n N_t^ie_i \in
\overline{\ka_+}$ (in particular it does not explode). Thus for
any $\alpha$, $(\int_0^t \coth \frac{(\alpha,X_s^W)}{2}ds,t\ge 0)$
is in fact a positive additive functional and its expectation is
therefore finite for each time $t\ge 0$ and for q.e. starting
point $x$ (q.e. stands for quasi everywhere, as explained in
\cite{FukOT}). But it results from \cite{CepLep2} Theorem $2.2$
that it is in fact true for all $x\in \overline{\ka_+}$. Indeed
there it is proved that $\E_x[\int_0^t|\nabla \log
\delta^{\frac{1}{2}}(X^W_s)|ds]<+\infty$. But $\nabla\log
\delta^{\frac{1}{2}}$ is equal to $k_\alpha \alpha \coth
\frac{\alpha}{2}+z$, where $z$ lies in the cone, let say $C^*$,
generated by the convex hull of $\kR^+$. Thus (since $-\alpha
\notin C^*$) there exists a constant $c>0$ (independent of $x$)
such that, $c\coth \frac{(\alpha,x)}{2}\le d(0,k_\alpha\alpha
\coth \frac{(\alpha,x)}{2}+C^*)\le |\nabla \log
\delta^{\frac{1}{2}}(x)|$, for $x\in \overline{\ka_+}$. Finally we
have proved the following result
\begin{prop}
\label{propositiondebase} The radial Heckman-Opdam process
$(X^W_t,t\ge 0)$ starting at $x\in \overline{\ka_+}$, is a
continuous semimartingale, and is the unique solution of the
following SDE
\begin{eqnarray}
\label{equationdebase} X^W_t=x+\beta_t+\sum_{\alpha \in
\kR^+}k_\alpha \frac{\alpha}{2} \int_0^t \coth
\frac{(\alpha,X_s^W)}{2}ds,\ t \ge 0, \end{eqnarray} where
$(\beta_t,t\ge 0)$ is a Brownian motion on $\ka$. Moreover for any
$t\ge 0$, any $x\in \overline{\ka_+}$ and any $\alpha \in \kR^+$,
$$\E_x \left[ \int_0^t \coth
\frac{(\alpha,X_s^W)}{2} ds\right] <+\infty .$$
\end{prop}
The uniqueness in law of the SDE (\ref{equationdebase}) is just a
consequence of the uniqueness of solutions to the martingale
problem associated to $(D,\kC(\ka)^W)$. In fact there is also
strong uniqueness. This results simply from the fact that $\coth$
is decreasing. Indeed if $(X,B)$ and $(X',B)$ are two solutions of
(\ref{equationdebase}), then for all $t\ge 0$,
$$\frac{d}{dt}(|X_t-X'_t|^2)=2(X_t-X'_t,\nabla\log
\delta^{\frac{1}{2}}(X_t)-\nabla\log
\delta^{\frac{1}{2}}(X'_t))\le 0,$$ which proves that $X$ and $X'$
are indistinguishable. With Theorem ($1.7$) p. $368$ in \cite{RY},
this implies also that each solution is strong.
\begin{rem} \label{remdebut} \emph{\begin{enumerate}
\item The finiteness of the
expectation in the proposition will be used in the next section
for the study of the jumps.
\item We could ask whether the processes considered by Cépa and Lepingle in
\cite{CepLep2} coincide with ours. In fact they prove existence of
a solution for the same EDS but with an additional local time
term. The question is therefore to know if this local time must be
$0$. Cépa and Lépingle have proved this for root systems of type
$A$. But we can prove it now for the other root systems. Indeed by
the Itô formula, their process is a solution of the martingale
problem associated to $D$, since for any $W$-invariant function
$f$, $(\nabla f(x),n)=0$ for all $x\in
\partial \ka_+$ and $n$ a normal vector. Thus they coincide with the radial HO-process whose local time on $\partial \ka_+$ is $0$.
\item In fact Proposition
\ref{propositiondebase} is also valid in the Dunkl setting (with
the same proof), where it was proved in the same time, but with a
completely different method, by Chybiryakov (see \cite{Chy2}).
\end{enumerate} }
\end{rem}
A first consequence of this proposition is an absolute continuity
relation between the laws of the radial HO-process and the
corresponding radial Dunkl process. More precisely, let $\pp^W$ be
the law of $(X^W_t,t\ge 0)$ with parameter $k$ on
$C(\R^+,\overline{\ka_+})$. Let
$\kR':=\{\frac{\sqrt{2}\alpha}{|\alpha|}\mid \alpha \in \kR\}$,
and if $\beta=\frac{\sqrt{2}\alpha}{|\alpha|} \in \kR'$, let
$k'_{\beta}:=k_\alpha+k_{2\alpha}$. Let $\Q^W$ be the law of the
radial Dunkl process $(Z^W_t,t\ge 0)$ associated to the root
system $\kR'$ and with parameter $k'$. Let $(L_t,t\ge 0)$ be the
process defined by
$$L_t:=\int_0^t \nabla
\log \frac{\delta^{\frac{1}{2}}}{\pi}(X_s)d\beta_s,\ t\ge 0,$$
where $(\beta_t,t\ge 0)$ is a Brownian motion under $\Q^W$, and
$$\pi(x)=\prod_{\beta \in \kR'}
(\beta,x)^{k'_\beta}.$$ As the function $x\mapsto
\frac{1}{x}-\coth(x)$ is bounded on $\R$ we get that, for all
$t\ge 0$, $\Q^W[\exp(\frac{1}{2}<L>_t)]<\infty$. Thus
$M:=\exp(L-\frac{1}{2}<L>)$, the stochastic exponential of $L$, is
a $\Q^W$-martingale. Moreover as mentioned in Remark
\ref{remdebut} we have also an explicit decomposition of the
radial Dunkl process as "Brownian motion + term with bounded
variation". Therefore, by using the Girsanov theorem \cite{RY}, we
get that for any $t\ge 0$, if $(\kF_t,t\ge 0)$ is the canonical
filtration on $C(\R^+,\overline{\ka_+})$, then
\begin{eqnarray}
\label{absocont} \pp^W_{|\kF_t}=M_t\cdot \Q^W_{|\kF_t}.
\end{eqnarray}
As a consequence we obtain for instance that when
$k_\alpha+k_{2\alpha} \ge 1/2$, then the HO-process starting at
any $x\in \ka_+$ a.s. does not touch the walls (i.e. the subspaces
of the type $\{\alpha =0\}$). This follows from the similar result
for the Dunkl processes proved in \cite{Chy}. Now if it starts at
some $x\in
\partial\ka_+$ then a.s., by the Markov property, it will never touch the walls in strictly positive
times (observe that $\mu(\partial \ka_+)=0$, thus at any $t>0$,
a.s. $X_t\in \ka_+$). \newline \newline We will now prove a law of
large numbers and a central limit theorem for the radial
Heckman-Opdam processes. These results are well known in the
setting of symmetric spaces of noncompact type (see for instance
Babillot \cite{Bab}).
\begin{prop}
\label{loigrdsnbres}The radial process satisfies the law of large
numbers
$$\lim_{t\to \infty} \frac{X^W_t}{t}\to \rho \text{  a.s.,}$$
and there is the convergence in $C(\R^+,\overline{\ka_+})$
$$(\frac{X^W_{tT}-\rho tT}{\sqrt{T}},t\ge 0) \xrightarrow[T\to
\infty]{} (\beta_t,t\ge 0).$$
\end{prop}
\textbf{Proof of the proposition:} The first step is to prove that
$(\alpha,X^W_t)\to +\infty$, when $t\to +\infty$, for all $\alpha
\in \kR^+$, or it is enough, for all the simple roots. We denote
by $\{\alpha_1,\dots,\alpha_n\}$ the set of simple roots. From
Proposition \ref{propositiondebase} we see that the radial process
(starting at $x$) satisfies for any $t\ge 0$,
\begin{eqnarray}
\label{EDSLGN} X^W_t= x + \beta_t+\rho t +\frac{1}{2}\sum_{\alpha
\in \kR^+} k_\alpha \alpha \int_0^t [\coth
\frac{(\alpha,X_s^W)}{2}-1] ds.
\end{eqnarray}
Let $u\in \overline{\ka_+}$. From (\ref{EDSLGN}) we get that for
all $t\ge 0$, $(u,X^W_t)-(u,x)-(u,\beta_t)\ge(u,\rho)t$, because
$\coth(x)\ge 1$ for $x\ge 0$. Thus $(u,X^W_t)\to +\infty$, when
$t\to +\infty$. In particular $(\rho,X^W_t)\to +\infty$. It
implies that $\max_{i=1,\dots,n}(\alpha_i,X^W_t)\to +\infty$. For
$t>0$, let $i_1,\dots,i_n$ be such that $(\alpha_{i_1},X^W_t)\ge
\dots\ge (\alpha_{i_n},X^W_t)$ (we forget the dependance in $t$ in
the notation). We prove now that $(\alpha_{i_2},X^W_t)\to
+\infty$. Let $\epsilon>0$ and let $T_0$ be such that
$\coth(\alpha_{i_1},X^W_t)-1\le \epsilon$ and $|x+\beta_t|\le
\epsilon t$ for $t\ge T_0$. Let $\kR_2$ be the root system
generated by $\{\alpha_{i_2},\dots,\alpha_{i_n}\}$ and let
$\rho_2=\sum_{\alpha \in \kR_2^+}k_\alpha \alpha$. Observe in
particular that if $\alpha \in \kR_2^+$, then $(\alpha,\rho_2)\ge
0$, whereas if $\alpha \notin \kR^+_2$, then
$\alpha-\alpha_{i_1}\in \kR^+$ and thus
$\coth\frac{(\alpha,X_t^W)}{2}\le 1+\epsilon$. Now from
(\ref{EDSLGN}) we get for $t\ge T_0$,
$$(\rho_2,X^W_t)\ge ((\rho,\rho_2)-\epsilon)t+f(t),$$
where $f(t)=\sum_{\alpha \in \kR^+} k_\alpha
(\alpha,\rho_2)\int_0^t[\coth\frac{(\alpha,X^W_s)}{2}-1]ds$. Hence
by our choice of $\rho_2$, we have for $t\ge T_0$, $f'(t)\ge
-C\epsilon$ for some constant $C>0$. Then we get another constant
$C'>0$ such that $f(t)\ge -C'-C\epsilon t$ for $t\ge T_0$. Thus we
conclude that $(\rho_2,X^W_t)\to+\infty$ and
$(\alpha_{i_2},X^W_t)\to +\infty$. In the same way we deduce that
$(\alpha_i,X^W_t)\to +\infty$ for all $1\le i\le n$. Eventually we
get immediately the law of large numbers from (\ref{EDSLGN}).
\newline For the second
claim of the proposition, we will show that a.s.
$$|\frac{1}{\sqrt{T}}(X^W_{tT}-x-\beta_{tT}-tT\rho)| \xrightarrow[T\to \infty]{} 0,$$ uniformly in $t\in \R^+$.
Let $\epsilon>0$. By the first claim, we know that there is some
$N$ such that for every $s\ge N$,
$|\coth\frac{(\alpha,X^W_s)}{2}-1|\le e^{-cs}$ for some strictly
positive constant $c$. Then we have
\begin{eqnarray*}
\frac{1}{\sqrt{T}}\int_0^{tT} [\coth \frac{(\alpha,X^W_s)}{2}-1]
ds &=& \frac{1}{\sqrt{T}}\int_0^{N\wedge tT}[\coth
\frac{(\alpha,X^W_s)}{2}-1] ds \\
&+& \frac{1}{\sqrt{T}}1_{(tT\ge N)} \int_N^{tT}[\coth
\frac{(\alpha,X^W_s)}{2}-1] ds. \end{eqnarray*} But the both
integrals can be made smaller than $\epsilon$ by choosing $T$
sufficiently large. The second claim follows using the scaling
property of the Brownian motion. \hfill $\square$

\section{Jumps of the process}
\label{sauts} We will now study the behavior of the jumps of the
Heckman-Opdam process. We use essentially the same tool as in
\cite{GY3} for the Dunkl processes, i.e. we use the predictable
compensators of some discontinuous functionals. However in our
setting we obtain a more precise result when
$k_\alpha+k_{2\alpha}\ge \frac{1}{2}$ for all $\alpha$. In fact in
this case, there is almost surely a finite random time, after
which the process does not jump anymore. This allows to prove for
such multiplicity $k$ a law of large numbers and a central limit
theorem for the HO-process.
\newline Let us first recall the definition of the Lévy kernel
$N(x,dy)$ of a homogeneous Markov process with a transition
semigroup $(P_t)_{t\ge0}$ and generator $\kD$ (see Meyer
\cite{Mey}). It is determined, for any $x\in \R^d$ by:
$$\kD f(x)=\lim_{t\to 0}\frac{P_tf(x)}{t}=\int_{\ka}N(x,dy)f(y),$$
for $f$ a function in the domain of the infinitesimal generator
which vanishes in a neighborhood of $x$. The following lemma
describes the Lévy kernel of the HO-process. It is an immediate
consequence of the explicit expression (\ref{explicit}) of $\kL$.
First we introduce some notation. If $I$ is a subset of $\kR^+$,
we denote by
$$\ka^I=\{x\in \ka \mid \forall \alpha \in I,\ (\alpha,x)=0 \},$$
the face associated to $I$. We denote also by $\kR_I$ the set of
positive roots which vanish on $\ka^I$.

\begin{lem}
The Lévy kernel of the HO-process has the following form
\begin{eqnarray*}
N(x,dy)= \left\{ \begin{array}{ccl}
                  \sum_{\alpha \in \kR^+}k_\alpha \frac{|\alpha|^2}{8}
                  \frac{\epsilon_{r_\alpha
                  x}(dy)}{\sinh^2\frac{(\alpha,x)}{2}} & if &
                  x\in \ka_{\text{reg}}\\
                 \sum_{\alpha \in \kR^+ \smallsetminus \kR_I }k_\alpha \frac{|\alpha|^2}{8}
                  \frac{\epsilon_{r_\alpha
                  x}(dy)}{\sinh^2\frac{(\alpha,x)}{2}} & if &
                  x\in \ka^I,\\
                  \end{array}
            \right.
\end{eqnarray*}
where $I$ is a subset of $\kR^+$, and for $x\in \ka$, $\epsilon_x$
is the Dirac measure in $x$.
\end{lem}
\begin{rem} \emph{The lemma implies that when there is a
jump at a random time $s$, i.e. $X_s \neq X_{s-}$, then almost
surely there exists $\alpha \in \kR^+$ such that $X_s=r_\alpha
X_{s-}$ (see \cite{GY3}). In this case we have
$$\Delta X_s:=X_s-X_{s-}=-(\alpha^\vee,X_{s-})\alpha.$$}
\end{rem}
Using the finiteness of the expectation of the time integrals
appearing in (\ref{equationdebase}), we can show that the sum over
any time interval of the amplitudes of the jumps is finite.
\begin{prop}
\label{ampliso} Let $(X_t,t\ge 0)$ be a Heckman-Opdam process. For
every $t>0$,
$$\E\left[\sum_{s\le t} |\Delta X_s|\right] < +\infty. $$
\end{prop}
\textbf{Proof of the proposition:} From the above remark we get
$$\sum_{s\le t} |\Delta X_s| = \sum_{\alpha \in \kR^+}\sum_{s\le t}
f_\alpha(X_{s-},X_s),$$ where
$$f_\alpha(x,y)=\frac{2}{|\alpha|}|(\alpha,x)|1_{(y=r_\alpha
x\neq x)}.$$ Now, the positive discontinuous functional
$\sum_{s\le t} f_\alpha(X_{s-},X_s)$ is compensated by the process
$\int_0^t ds\int_{\ka}N(X_{s-},dy)f_\alpha(X_{s-},y)$. As a
consequence, the proposition will be proved if we know that the
expectation of the compensator is finite at all time $t\ge 0$.
Thus we have to show that for every $\alpha \in \kR^+$,
$$\E
\left[\int_0^t\left| \frac{(\alpha,X_s)}{\sinh^2
\frac{(\alpha,X_s)}{2}}\right|ds \right] <+\infty.$$ But for every
$x>0$, $\frac{x}{\sinh^2x} \le 2\coth x$. Therefore the above
condition follows from Proposition \ref{propositiondebase}. \hfill
$\square$
\newline
\newline
For $\alpha \in \kR^+$, we denote by $(M^\alpha_t,t\ge 0)$ the
process defined for $t\ge 0$ by:
\begin{eqnarray}
\label{malpha} M^\alpha_t= \sum_{s\le t}
-(\alpha^\vee,X_{s-})1_{(r_\alpha X_{s-}=X_s)}+\frac{k_\alpha}{4}
\int_0^t \frac{(\alpha,X_s)}{\sinh^2\frac{(\alpha,X_s)}{2}}ds.
\end{eqnarray}
By the martingale characterization of $(X_t,t\ge 0)$, we know that
$(f(X_t),t\ge 0)$ is a local semimartingale for all $f\in
C_c^\infty(\ka)$. Thus $(X_t,t\ge 0)$ itself is a local
semimartingale. In the next proposition we give its explicit
decomposition. It is the analogue of a result of Gallardo and Yor
on the decomposition of the Dunkl processes. The proof is very
similar (and uses Proposition \ref{ampliso}), so we refer to
\cite{GY3} for more details.
\begin{prop}
\label{decompositionsm} We have the following semimartingale
decomposition:
$$X_t=\beta_t+\sum_{\alpha \in \kR^+}M_t^\alpha \alpha+A_t,$$
for $t\ge 0$, where $$A_t=\sum_{\alpha \in
\kR^+}\frac{k_\alpha}{2} \alpha \int_0^t \left[\coth
\frac{(\alpha,X_s)}{2}-
\frac{(\alpha,X_s)}{2\sinh^2\frac{(\alpha,X_s)}{2}}\right]ds,$$
and the $M^\alpha$'s are purely discontinuous martingales given by
(\ref{malpha}) which satisfy
$$[M^\alpha,M^\beta]_t=0, \text{ if } \alpha \neq \beta,$$ and
$$<M^\alpha>_t=\frac{k_\alpha}{4|\alpha|^2}\int_0^t\frac{(\alpha,X_s)^2}{\sinh^2\frac{(\alpha,X_s)}{2}}ds.$$
\end{prop}

Another interesting property of the jumps is that, when $k_\alpha
+ k_{2 \alpha}\ge 1/2$ for all $\alpha$, and when the starting
point lies in $\ka_{\text{reg}}$, in which case the HO-process
does not touch the walls, the number of jumps $N_t$ up to a fixed
time $t$ is a.s. finite. Indeed otherwise the paths of the
trajectories would not be càdlàg. Therefore the sequence of
stopping times $T_n=\inf \{t>0, N_t\ge n\}$ converges a.s. to
$+\infty$ when $n$ tends to infinity. Thus $(N_t,t\ge 0)$ is a
locally integrable (because locally finite) increasing process. We
will deduce from this observation and a general result of
\cite{Jac} a more precise result. For $t\ge 0$, we denote by $w_t$
the element of $W$ such that $X_t=w_tX_t^W$.
\begin{prop}
Assume that $k_\alpha + k_{2\alpha} \ge 1/2$ for all $\alpha$.
\begin{enumerate}
\item If the starting point lies in $\ka_{\text{reg}}$, then a.s.
$$\sup_{t\ge 0} N_t < +\infty.$$
\item For any starting point in $\ka$, $w_t$ converges a.s. to $w_\infty\in W$ when $t\to \infty$. If the process starts from zero, then the law of $w_\infty$ is the uniform probability on $W$.
\item When $T\to \infty$, the sequences $(\frac{1}{T}X_{tT},t\ge 0)$ and $(\frac{1}{\sqrt{T}}(X_{tT}-w_{tT} \rho tT),t\ge 0)$ converge in
law in $\D(\R^+,\ka)$ respectively to $(w_\infty\rho t,t\ge 0)$,
and to a Brownian motion $(\beta_t,t\ge 0)$.
\end{enumerate}
\end{prop}
\textbf{Proof of the proposition:} Let us begin with the first
claim. As in the preceding proposition we use the following result
of Meyer about the Lévy kernel: the positive discontinuous
functional $(N_t,t\ge 0)$ can be compensated by the predictable
process $(\sum_{\alpha \in \kR^+} k_\alpha \frac{|\alpha|^2}{8}
\int_0^t \frac{1}{\sinh^2\frac{(\alpha,X_s)}{2}}ds,t\ge 0)$. Now
from the law of large numbers (Proposition \ref{loigrdsnbres}) we
deduce that this compensator converges a.s. to a finite value when
$t\to \infty$. Thus the corollary $(5.20)$ p.$168$ of \cite{Jac}
gives the result. Then the second point is simply a consequence of
the first point and of the Markov property. The assertion that the
limit law is uniform when the process starts from zero results
from the fact that for any $w\in W$, $\kD$ remains unchanged if we
replace $\kR^+$ by $w\kR^+$. The first convergence result of the
last point is straightforward with the second point and
Proposition \ref{loigrdsnbres}. For the second convergence result,
we can use Proposition \ref{decompositionsm}. Indeed it says that
for all $t>0$ and $T>0$, \begin{eqnarray*} \frac{X_{tT}-w_{tT}\rho
tT}{\sqrt{T}}&=&
\frac{\beta_{tT}}{\sqrt{T}}+\frac{1}{\sqrt{T}}\sum_{s\le tT}\Delta
X_s\\
&+&\sum_{\alpha \in \kR^+}k_\alpha
\frac{\alpha}{2\sqrt{T}}\int_0^{tT} \left[\coth
\frac{(\alpha,X_s)}{2}-\epsilon_{tT}^\alpha\right]ds,
\end{eqnarray*}
where the $\epsilon_{tT}^\alpha\in \{\pm 1\}$ are defined by
$w_{tT}\rho = \sum_{\alpha \in \kR^+}k_\alpha \epsilon_{tT}^\alpha
\alpha$, or equivalently by $\epsilon_{tT}^\alpha \alpha \in
w_{tT}\kR^+$. But by the second point of the proposition, we know
that a.s. there exists a random time after which the process stays
in the same chamber, which is $w_\infty \ka_+$. Moreover for all
$s>0$ and all $\alpha \in \kR^+$,
$\epsilon_s^\alpha(\alpha,X_s)\ge 0$. Thus by Proposition
\ref{loigrdsnbres}
$\coth\frac{(\alpha,X_s)}{2}-\epsilon_{tT}^\alpha$ tends to $0$
exponentially fast when $s\to \infty$ (and $s\le tT$). Then a.s.
$\sum_{\alpha \in \kR^+}k_\alpha
\frac{\alpha}{2\sqrt{T}}\int_0^{tT} \left[\coth
\frac{(\alpha,X_s)}{2}-\epsilon_{tT}^\alpha\right]ds$ tends to $0$
when $T\to \infty$, uniformly in $t\in \R^+$. In the same way, by
Proposition \ref{ampliso}, a.s. for any $A>0$, $\sum_{s\le
A}|\Delta X_s|<+\infty$. By the second point we know that a.s.
after some time there is no more jumps, thus a.s.
$\frac{1}{\sqrt{T}}\sum_{s\le tT}|\Delta X_s|$ tends to $0$ when
$T\to \infty$, uniformly in $t\in \R^+$. This proves the desired
result by the scaling property of the Brownian motion. \hfill
$\square$

\section{Convergence to the Dunkl processes} In this section we will show that when it is well normalized, the
HO-process of parameter $k>0$ converges to a certain Dunkl process
$(Z_t,t\ge 0)$. The proof uses a general criteria for a sequence
of Feller processes with jumps, which can be found in \cite{EK}
for instance. Roughly speaking it states that it just suffices to
prove the convergence of the generator of these processes on a
core of the limit. Let us notice that the convergence of the
normalized radial HO-process to the corresponding radial Dunkl
process is more elementary. It could be proved essentially by
using that the laws of the radial HO-process and the radial Dunkl
process are absolutely continuous, and that the normalized
Radon-Nikodym derivative tends to $1$. Let us also observe that
the convergence of the normalized radial process has a natural
geometric interpretation in the setting of symmetric spaces of
noncompact type. Indeed in this setting the radial Dunkl process
is just the radial part of the Brownian motion on the tangent
space (or the Cartan motion group, see the more precise
description by De Jeu \cite{J2}, and in \cite{ABJ} in the complex
case). From the analytic point of view, it also illustrates the
more conceptual principle, that the Dunkl (also called rational)
theory is the limit of the Heckman-Opdam (or trigonometric)
theory, when "the curvature goes to zero".
\newline We denote by $(X^T_t,t\ge 0)$ the normalized HO-process, which is defined for $t\ge 0$ and $T>0$ by:
$$X^T_t=\sqrt{T}X_{\frac{t}{T}}.$$
We recall that $\kR'=\{\frac{\sqrt{2}\alpha}{|\alpha|},\ \alpha
\in \kR\}$, and that for $\beta= \frac{\sqrt{2}\alpha}{|\alpha|}
\in \kR'$, $k'_{\beta}=k_\alpha+k_{2\alpha}$.
\begin{theo}
\label{convergencedunkl} When $T\to \infty$, the normalized
HO-process $(X^T_t,t\ge 0)$ with parameter $k$ starting at $0$
converges in distribution in $\D(\R^+,\ka)$ to the Dunkl process
$(Z_t,t\ge 0)$ associated with $\kR'$ and with parameter $k'$
starting at $0$.
\end{theo}
\textbf{Proof of the theorem:} First it is well known that the
process $(X^T_t,t\ge 0)$ is also a Feller process with generator
$\kL^T$ defined for $f\in C^2(\ka)$ and $x\in \ka_{\text{reg}}$,
by $\kL^Tf(x)= \frac{1}{T}(\kL g)(\frac{x}{\sqrt{T}})$, where
$g(x)=f(\sqrt{T}x)$. Thus a core of $\kL^T$ is, like for $\kL$ and
$\kL'$, the space $C_c^\infty(\ka)$. Moreover it is
straightforward that for any $f\in C_c^\infty(\ka)$, $\kL^Tf$
converges uniformly on $\ka$ to $\kL'f$. Therefore we can apply
Theorem $6.1$ p.$28$ and Theorem $2.5$ p. $167$ in \cite{EK},
which give the desired result. \hfill $\square$

\section{The $F_0$-process and its asymptotic behavior} We introduce here and study a generalization
of the radial part of the Infinite Brownian Loop (abbreviated as
I.B.L.) introduced in \cite{ABJ}. Let
$\tilde{F}_0(x,t):=F_0(x)e^{\frac{|\rho|^2}{2}t}$, for $(x,t)\in
\ka \times [0,+\infty)$. Then $\tilde{F}_0$ is harmonic for the
operator $\partial_t+\kD$ which is the generator of $(X_t,t)_{t\ge
0}$. We define now the processes $(Y_t,t)_{t\ge 0}$ as the
relativized $\tilde{F}_0$-processes in the sense of Doob of
$(X_t,t)_{t\ge 0}$. By abuse of language we will call $(Y_t,t\ge
0)$ the $F_0$-process. We denote by $(Y_t^W,t\ge 0)$ its radial
part, that we will call the radial $F_0$-process. For particular
values of $k$ it coincides with the radial part (in the Lie group
terminology) of the I.B.L. on a symmetric space.
\newline
\newline
The goal of this section is to generalize some results of
\cite{ABJ} and \cite{BJ}, for any $k>0$. Essentially we first
prove the convergence of the HO-bridge of length $T$ around $0$,
i.e. the HO-process conditioned to be equal to $0$ at time $T$, to
the $F_0$-process starting at $0$, when $T$ tends to infinity.
Then we prove the convergence of the normalized $F_0$-process to a
process whose radial part is the intrinsic Brownian motion, but
which propagates in a random chamber (independently and uniformly
chosen). We begin by the following lemma:
\begin{lem} Let $x,a \in \overline{\ka_+}$. When $T\to \infty$,
$$\frac{p^W_{T-t}(x,a)}{p^W_T(a,a)}\to
\frac{F_0(x)}{F_0(a)}e^{\frac{t}{2}|\rho|^2}.$$
\end{lem}
\textbf{Proof of the lemma:} We need the integral formula of the
heat kernel:
$$p^W_t(x,y)=\int_{i\ka}e^{-\frac{t}{2}(|\la|^2+|\rho|^2)}F_\la(x)F_{-\la}(y)d\nu'(\la),\ x,y\in \ka_+.$$
We make the change of variables $u:=\la (T-t)$ for $p^W_{T-t}$,
and $v:=\la T$ for $p^W_T$. Then we let $T$ tend to $+\infty$ and
the result follows. \hfill $\square$

\begin{prop} \label{pontho} Let $(X_t^{0,T},t\ge 0)$ be the HO-bridge around
$0$ of length $T$. Then when $T\to +\infty$, it converges in
distribution in $\D(\R^+,\ka)$ to the
$F_0$-process starting at $0$.
\end{prop}
\textbf{Proof of the proposition:} We know that for any $t\ge 0$,
and any bounded $\kF_t$-measurable function $F$,
\begin{eqnarray*}
\E \left[ F(X^{0,T}_s,0\le s\le t)\right] &=& \E \left[ F(X_s,0\le s\le
t)\frac{p_{T-t}(X_t,0)}{p_T(0,0)} \right]\\
&=& \E \left[ F(X_s,0\le s\le
t)\frac{p^W_{T-t}(X^W_t,0)}{p^W_T(0,0)} \right]
\end{eqnarray*} The second equality results from the fact that $p_t(x,0)=\frac{1}{|W|}p_t^W(x,0)$, for all $x\in \ka$ and all $t\ge 0$ (see \cite{Schap}). Moreover since
$F_\la$ is bounded (cf \cite{O}) and the measure $\nu'$ is
positive, we see from the integral formula of $p_t^W$, that there
exists a constant $C$ such that, $p_t^W(x,y)\le C p^W_t(0,0)$, for all $t>0$ and all $x,y\in
\overline{\ka_+}$. It follows
that $\frac{p^W_{T-t}(x,a)}{p^W_T(a,a)}$ is a bounded
function of $(x,T)\in \overline{\ka_+} \times [1,\infty)$. Then we
get the result from the preceding lemma and the dominated
convergence theorem. \hfill $\square$
\begin{rem}
\emph{The same
proof shows in fact that for all $a\in \overline{\ka_+}$, the radial HO-bridge of length $T$ around $a$
converges in law to the radial $F_0$-process starting from $a$.}
\end{rem}
The next result is an important technical lemma:
\begin{lem}
\label{eulerfo} There exist two Bessel processes $(R_t,t\ge 0)$
and $(R'_t,t\ge 0)$ (not necessarily with the same dimension),
such that a.s. $|R_0|=|R'_0|=|Y_0|$ and for any $t \ge 0$,
$$R^2_t\le |Y_t|^2 \le R'^2_t.$$
\end{lem}
\textbf{Proof of the lemma:} First the $F_0$ process and its
radial part have the same norm, hence it suffices to prove the
result for $(Y^W_t,t\ge 0)$. Next we can follow exactly the same proof
as in \cite{ABJ}. We recall it for the convenience of the reader.
We know that $(Y^W_t,t\ge 0)$ is solution of the SDE
$$Y^W_t=Y^W_0+\beta_t+\int_0^t\nabla \log (\delta^{\frac{1}{2}}F_0)(Y^W_s) ds,$$
where $(\beta_t,t\ge 0)$ is a Brownian motion. By Itô formula we
get $$|Y^W_t|^2=|Y^W_0|^2+2\int_0^t(Y^W_s,d\beta_s)+tn+2\int_0^t
E[\log(\delta^{\frac{1}{2}}F_0)](Y^W_s)ds,$$ where $E=\sum_{i=1}^n
x_i\frac{\partial}{\partial x_i}$ is the Euler operator on $\ka$.
And it was shown in \cite{Schap} that
$E[\log(\delta^{\frac{1}{2}}F_0)]$ is positive and bounded on
$\ka$. Thus we can conclude by using comparison theorems. \hfill
$\square$

\begin{cor} Let $(Y_t,t\ge 0)$ be the $F_0$-process. Then almost
surely,
$$\lim_{t\to \infty} \frac{|Y_t|}{t}=0.$$
More precisely (law of the iterated logarithm), a.s.
$$\limsup_{t\to \infty}\frac{|Y_t|}{\sqrt{2t\log \log t}}=1.$$
\end{cor}
\textbf{Proof of the corollary:} It follows from the preceding
lemma and known properties of the Bessel processes. \hfill
$\square$
\newline
\newline
In the complex case, i.e. when $2k$ is equal to the multiplicity
function on some complex Riemannian symmetric space of noncompact
type (or equivalently when $\kR$ is reduced and $k=1$), then it
was proved in \cite{ABJ} that the radial $F_0$-process coincides with the
intrinsic Brownian motion. It was also proved in \cite{ABJ} that
in the real case, i.e. for other choices of $k$, a
normalization of the radial $F_0$-process converges to the intrinsic
Brownian motion. The next theorem gives a generalization of this result for any
multiplicity $k>0$ and for the (non radial) $F_0$-process also. We denote by $(Y^{T}_t,t\ge 0)$ the process
defined for $t\ge 0$ and $T>0$, by
$$Y^{T}_t:=\frac{1}{\sqrt{T}}Y_{tT},$$
and we denote by $(Y^{W,T}_t,t\ge 0)$ its radial part.
Let $(I_t,t\ge 0)$ be the intrinsic Brownian motion starting from $0$. We denote by $(I^*_t,t\ge
0)$ the continuous process starting from $0$, whose radial part is
the intrinsic Brownian motion, but which propagates in a random
chamber $w \ka_+$, where $w$ is chosen independently and with
respect to the uniform probability on $W$. This is a typical
example of a strong Markov process which is not Feller (for
instance it does not satisfy the Blumenthal's zero-one law). We
have
\begin{theo}
\label{convergencefo} Let $k>0$. The normalized $F_0$-process
$(Y^T_t,t\ge 0)$ starting at any $x\in \ka$ converges in
distribution in $\D(\R^+,\ka)$ to $(I^*_t,t\ge 0)$.
\end{theo}
\textbf{Proof of the theorem:} The first step is to prove that $(Y^{W,T}_t,t\ge 0)$ converges
in law in $C(\R^+,\overline{\ka_+})$ to $(I_t,t\ge 0)$. Thanks to Lemma \ref{eulerfo} we
can use the same proof as for Theorem $5.5$ in \cite{ABJ}. The
result about the convergence of the semigroup needed in the proof
was established in \cite{Schap}. Now let $\pp^0$ be the law of the
$F_0$-process, and let $\pp$ be the one of the HO-process. By
definition we have the absolute continuity relation
$$\pp^0_{|\kF_t}=\frac{F_0(Y_t)}{F_0(x)}e^{\frac{|\rho|^2}{2}t}\cdot
\pp_{|\kF_t}.$$ Since $F_0$ is bounded (cf \cite{O}) Proposition
\ref{ampliso} implies that for any $t>0$, $\E^0[\sum_{s\le t}
|\Delta Y_s|]<+\infty$. Thus by Girsanov theorem (see \cite{RY})
and Proposition \ref{decompositionsm} we get the semimartingale
decomposition of the $F_0$-process:
\begin{eqnarray*}
Y_t &=& x+\beta_t+\underbrace{\sum_{\alpha \in \kR^+} M^\alpha_t
\alpha}_{:=M_t} +
\int_0^t\nabla \log \delta^{\frac{1}{2}}F_0(Y_s)ds\\
&-& \sum_{\alpha \in \kR^+} \frac{k_\alpha}{2} \alpha \int_0^t
\frac{(\alpha,Y_s)}{2\sinh^2\frac{(\alpha,Y_s)}{2}}d s,
\end{eqnarray*}
where $(\beta_t,t\ge 0)$ is a $\pp^0$-Brownian motion and the $M^\alpha$'s are defined by (\ref{malpha}) (with $Y_s$
in place of $X_s$). We set $M^T_t:=\frac{1}{\sqrt{T}}M_{tT}$. By
Proposition \ref{decompositionsm} we know that
$$<M^T>_t:=\sum_{i=1}^n <M^T_i>_t=\sum_{\alpha \in
\kR^+}\frac{k_\alpha}{4}\int_0^t
\frac{(\alpha,\sqrt{T}Y^T_s)^2}{\sinh^2 \frac{\sqrt{T}}{2}
(\alpha,Y^T_s)}ds,$$ where $M^T_i$ is the $i^{\text{th}}$
coordinate of $M^T$ in the canonical basis. Now
for all $w\in W$ the preceding sum remains unchanged if $\kR^+$ is replaced by $w\kR^+$. Therefore we can replace $Y^T_s$ by $Y^{W,T}_s$ in the last equality.
Moreover $\sinh x\ge x + \frac{x^3}{6}$ on $\R^+$. Hence
$$<M^T>_t\ \le \sum_{\alpha \in
\kR^+}k_\alpha\int_0^t\frac{1}{1+\frac{T}{12}(\alpha,Y^{W,T}_s)^2}ds.$$
Thus by using the first step, we see that for any fixed $t>0$,
$\E[<M^T>_t]\to 0$ when $T\to +\infty$. It implies by Doob's
$L^2$-inequality (see \cite{RY}), that $(M^T_t,t\ge 0)$ converges
in law in $\D(\R^+,\ka)$ to $0$. Now the triangular inequality
implies that for any $A>0$, $\epsilon>0$ and $\alpha>0$,
\begin{eqnarray*}
\pp\left[\sup_{|s-t|\le \epsilon,\ s\le t\le A}|Y_s^T-Y_t^T|\ge
\alpha\right] &\le & \pp\left[\sup_{|s-t|\le \epsilon,\ s\le t\le
A}|Y_s^{W,T}-Y_t^{W,T}|\ge
\frac{\alpha}{2N}\right]\\
&+& \pp\left[\sum_{s\le A}|\Delta Y_s^T|\ge
\frac{\alpha}{2}\right],
\end{eqnarray*} where $N$ is the number of Weyl chambers. Thus using that $\Delta Y^T=\Delta M^T$, tightness of $(Y^{W,T}_t,t\ge 0)$, and standard results (see
Theorem $3.21$ p.$314$ and Proposition $3.26$ p.$315$ in
\cite{JacShir} for instance), we see that the sequence
$(Y^T_t,t\ge 0)$ is C-tight in $\D(\R^+,\ka)$, i.e. it is tight
and any limit law of a subsequence is supported on $C(\R^+,\ka)$.
By the first step each limit has a radial part equal to the
intrinsic Brownian motion. Since we know that the intrinsic
Brownian motion does not touch the walls (in strictly positive
time), each limit process is necessarily of the type $(w
I_t,t\ge 0)$ where $w$ is some random variable on $W$. Thus in
order to identify the limit, we need to prove that the law of
$w$ has to be the uniform probability on $W$, and that $w$ is
independent of the radial part. For the law of $w$
first, let us just observe that when the process starts from $0$,
the result is immediate since by $W$-invariance of $\kD$ and
$F_0$, the law of the $F_0$-process is $W$-invariant, and thus the
law of any limit also. However when $x\neq 0$ we can not argue
like this, and we need to prove for instance that the law of
$Y^T_1$ converges to the law of $I_1^*$, when $T\to \infty$. Let
$f:\ka\to \R$ be continuous and bounded. We have
\begin{eqnarray*}
\E_{\frac{x}{\sqrt{T}}}\left[f(\frac{Y_T}{\sqrt{T}})\right]&=&
\int_\ka
p_T(\frac{x}{\sqrt{T}},\sqrt{T}y)\frac{F_0(\sqrt{T}y)}{F_0(\frac{x}{\sqrt{T}})}e^{\frac{|\rho|^2}{2}T}f(y)d\mu(\sqrt{T}y).
\end{eqnarray*}
Then it results from the asymptotic of
$p_T(\frac{x}{\sqrt{T}},\sqrt{T}y)$ and of $F_0(\sqrt{T}y)$ proved
in \cite{Schap}, that
$$p_T(\frac{x}{\sqrt{T}},\sqrt{T}y)\frac{F_0(\sqrt{T}y)}{F_0(\frac{x}{\sqrt{T}})}e^{\frac{|\rho|^2}{2}T}d\mu(\sqrt{T}y)
\to \text{const}\cdot e^{-\frac{|y|^2}{2}}\prod_{\alpha\in
\kR^+}(\alpha,y)^2dy,$$ which gives the density of the law of
$I^*_1$ (see \cite{ABJ} for instance for the law of $I_1$). We
conclude by Sheffé's lemma. Now the only missing part is the
independence of $w$ and $(I_t,t\ge 0)$. Observe first that since
any limit process is adapted, $w$ is $\kF_{0^+}$-measurable. On
the space $\D(\R^+,\ka)$, we denote by $(\kF_t^W)_t\ge 0$ the
natural filtration of the radial process. We know that $(I_t,t\ge
0)$ is an $(\kF^W_t)_{t\ge 0}$-Markov process. We will prove that
it is also an $(\kF_t)_{t\ge 0}$-Markov process. Indeed since it
is a.s. continuous and equal to $0$ at time $0$, this will imply the
independence with $w$. We know that $Y^{W,T}$ is an
$(\kF_t)_{t\ge 0}$-Markov process, since it is the projection of
$Y^T$. We denote by $P^{W,T}$ the semigroup of $Y^{W,T}$, and by
$Q^W$ the semigroup of $(I_t,t\ge 0)$. In particular we already
know that for $t\ge 0$ and any continuous and bounded function $f$,
$P_t^{W,T}f$ converges simply to $Q_t^Wf$ when $T\to \infty$. For
$s<t$, and $f$ and $g$ continuous and bounded functions, we have
$\E[f(Y^{W,T}_t)g(Y^T_s)]=\E[P^{W,T}_{t-s}f(Y^{W,T}_s)g(Y^T_s)]$.
For a suitable subsequence of $T$, the first term tends to
$\E[f(I_t)g(\nu I_s)]$, and the second term tends to
$\E[Q^W_{t-s}f(I_s)g(\nu I_s)]$, which implies the desired result.
This finishes the proof of the theorem. \hfill $\square$
\newline
\newline
We can now prove a generalization of a result of Bougerol and
Jeulin \cite{BJ}. Let $(R^{0,T}_t,0\le t\le 1)$ be the normalized
HO-bridge of length $T$ around $0$. It
is defined for $t\ge 0$ by
$$R^{0,T}_t=\frac{1}{\sqrt{T}}X^{0,T}_t.$$
\begin{theo}
\label{cras} Let $k>0$. When $T\to \infty$, the process
$(R^{0,T}_t,0\le t \le 1)$ converges in distribution in
$\D(\R^+,\ka)$ to the bridge $(I_t^{\{*,0,0,1\}},0\le
t\le 1)$ of length $1$ associated to $(I^*_t,0\le t\le 1)$.
\end{theo}
\textbf{Proof of the theorem:} Here again we can follow the same
proof as in \cite{BJ}. We just need to take care that the estimate
of the heat kernel in Proposition $5.3$ in \cite{Schap}, is
uniform when $y$ lies in some compact of $\ka_+$. But this results
directly from the proof of this proposition. \hfill $\square$

\begin{rem}\emph{We can define similarly the normalized radial HO-bridge around any $a\in \overline{\ka_+}$. With the same proof, we can also prove that it converges to the bridge $(I_t^{\{0,0,1\}},0\le t \le 1)$ of length $1$ associated to the intrinsic Brownian motion starting from $0$. Let us just notice that in dimension $1$ this is the
bridge of a Bessel-$3$, which is also the normalized Brownian
excursion. Thus we do recover the result of \cite{BJ2}}.
\end{rem}

\noindent \textbf{Acknowledgments:} I very much thank my advisors
Jean-Philippe Anker and Philippe Bougerol for their help and many
suggestions. I wish to thank also Marc Yor for his comments of a
previous version, and Oleksandr Chybiryakov and Emmanuel Cépa for 
fruitful discussions.

\end{document}